\documentclass{article}%
\usepackage{amsmath}
\usepackage{amsfonts,amsthm}
\usepackage{amssymb}
\usepackage{graphicx}
\usepackage{algorithm}
\usepackage{algorithmic}
\usepackage{fancybox}
\usepackage{xcolor}
\usepackage{tikz}
\usepackage{float}%

\providecommand{\U}[1]{\protect \rule{.1in}{.1in}}
\newtheorem{theorem}{Theorem}
\newtheorem{corollary}[theorem]{Corollary}
\newtheorem{lemma}[theorem]{Lemma}
\newtheorem{proposition}[theorem]{Proposition}
\theoremstyle{definition}
\newtheorem{remark}[theorem]{Remark}
\newtheorem{example}[theorem]{Example}
\def\dj{d\kern-0.4em\char"16\kern-0.1em}

\DeclareMathOperator{\main}{main}

\evensidemargin=\oddsidemargin
\begin{document}

\title{The main vertices of a star set and related graph parameters \footnote{To appear in Discrete Mathematics.}}
%The main vertices of a star set - a simplex like approach
\date{}
\author{Milica An\dj eli\'c \thanks{Department of Mathematics, Kuwait University, Safat 13060 Kuwait.
E-mail: milica.andelic@ku.edu.kw ORCID: http://orcid.org/0000-0002-3348-1141.}
\and Domingos M. Cardoso \thanks{Center for Research and Development in Mathematics
and Applications, Department of Mathematics, University of Aveiro, Campus de Santiago, 3810-193 Aveiro, Portugal.
E-mail: dcardoso@ua.pt ORCID: http://orcid.org/0000-0001-6239-3557}
\and Slobodan K.~Simi\' c \thanks{We started this research in 2015 with Slobodan K.~Simi\' c from Mathematical Institute of SANU, Belgrade, Serbia, who passed away in May 2019.}
\and Zoran Stani\'c \thanks{Faculty of Mathematics, University of Belgrade, Studentski trg 16, 11 000 Belgrade, Serbia.
E-mail: zstanic@matf.bg.ac.rs ORCID: http://orcid.org/0000-0002-4949-4203.}}

\maketitle

\begin{abstract}
A vertex $v \in V(G)$ is called $\lambda$-main if it belongs to a star set $X \subset V(G)$ of the
eigenvalue $\lambda$  of a graph $G$ and this eigenvalue is main for the graph obtained from $G$
by deleting all the vertices in $X \setminus \{v\}$; otherwise, $v$ is $\lambda$-non-main.
Some results concerning main and non-main vertices of an eigenvalue  are deduced. For a main
eigenvalue $\lambda$ of a graph $G$, we introduce the minimum and maximum number of $\lambda$-main
vertices in some $\lambda$-star set of $G$ as new graph invariant parameters. The determination of
these parameters is formulated as a combinatorial optimization problem based on a simplex-like approach.
Using these and some related parameters we develop new spectral tools that can be used in the research
of the isomorphism problem. Examples of graphs for which the maximum number of $\lambda$-main vertices
coincides with the cardinality of a $\lambda$-star set are provided.

\medskip


\smallskip


\end{abstract}

\section{Introduction}

Throughout this paper we consider undirected simple graphs $G$ with vertex set $V(G)=\{1,2,\ldots, n\}$
and edge set $E(G)$. An edge linking the vertices $i$ and $j$ of $V(G)$ is denoted by $ij\in E(G)$, and in
this case we say that $i$ and $j$ are \textit{adjacent}. For each vertex $i \in V(G),$ $N_{G}(i)$
denotes its \textit{neighbourhood}, that is the set of vertices of $G$ which are adjacent
to~$i$ and $\left \vert N_{G}(i)\right \vert$ is called the \textit{degree} of  $i$ and denoted by $d_G(i)$.
Given $S \subseteq V(G)$, the subgraph of $G$ \textit{induced} by $S$ is denoted by $G[S]$ and is such that
$V(G[S])=S$ and $E(G[S])=\{ij\in E:i,j\in S\}$.

The \textit{adjacency} matrix $A_{G}=[a_{ij}]$ of $G$ is the symmetric matrix such that $a_{ij}=1$ if $ij\in E(G)$
and $0$, otherwise.  The multiset of eigenvalues of $A_{G}$ (called
the \textit{spectrum} of~$G$) is defined as $\sigma(G) = \big\{\mu_1^{[k_1]}, \mu_2^{[k_2]}, \ldots, \mu_m^{[k_m]}\big\}$, where $\mu_i^{[k_i]}$ means that the eigenvalue $\mu_i$ appears repeated $k_i$ times in the spectrum of $G$.
The \textit{eigenspace} of  $\lambda \in \sigma(G)$ is denoted by $\mathcal{E}_G(\lambda)$, that is,
$\mathcal{E}_G(\lambda)=\ker(A_{G}-\lambda I_{n})$, where $I_{n}$ is the $n\times n$ identity matrix,
considering a square matrix $M$, $\ker(M)$ is the \textit{kernel} (or null space) of $M$.

Each of the  eigenvalues $\mu_{1}, \mu_{2}, \ldots, \mu_{m}$ of a graph $G$ whose
eigenspace $\mathcal{E}_G(\mu _i)$ is not orthogonal to the all-1 vector with $n$ entries $\textbf{j}_n$ is
said to be \textit{main}; otherwise, it is \textit{non-main}. The concept of main (non-main)
eigenvalue was introduced by Cvetkovi\' c in \cite{cvetkov70} and further investigated in several publications.
%{\color{red}Among many published results, it is useful to highlight the following proposition.
%\begin{proposition}\cite{eig}
%Let $\mu_i$ adn $\mu_j$ be algebraic conjugate eigenvalues of a graph $G$.
%Then $\mu_i$ is a main eigenvalue of $G$ if and only if $\mu_j$ is a main eigenvalue of $G$.
%\end{proposition}}
A survey on main eigenvalues is exposed by Rowlinson in \cite{rowmain}.

The remaining part of the paper is organized as follows. In Section~\ref{sec_2} we give some preliminary results.
In Section~\ref{sec_3} the concepts of main and non-main vertices are introduced and several theoretical results are established. In particular, it is proved that, for some main eigenvalue $\lambda$, a particular vertex can be
$\lambda$-main for some star set and $\lambda$-non-main for another star set. In Section~\ref{sec_4} the graph invariants related to the maximum and the minimum number of $\lambda$-main vertices are introduced and their determination is formulated as a combinatorial optimization problem based on a simplex-like approach. Furthermore, these invariants are related to the graph isomorphism problem. In Section~\ref{sec:max} we construct some examples of graphs in which all vertices of a fixed $\lambda$-star set are $\lambda$-main. Some open problems we observed during the research are selected in Section~\ref{sec:problems}. A computation that supports some results of Section~\ref{sec_4} is separated in the Appendix.

\section{Preliminary results on star sets and star complements}\label{sec_2}
We first recall some basic concepts of the theory of star sets. For more details we refer to \cite[pp. 136--141]{CvetRowSim2010}.

Considering a graph $G$ with $n$ vertices and an eigenvalue $\lambda \in \sigma(G),$ let $P$ be the matrix of the orthogonal projection of $\mathbb{R}^{n}$ onto $\mathcal{E}_G(\lambda)$ with respect to the standard orthonormal basis $\{\mathbf{e}_{1}, \mathbf{e}_2, \ldots, \mathbf{e}_{n}\}$ of $\mathbb{R}^{n}$. Then the set of vectors $P\mathbf{e}_{j}$ $(1\leq j\leq n)$ spans $\mathcal{E}_G(\lambda)$, and therefore there exists $X \subseteq V(G)$ such that the vectors $P\mathbf{e}_{j}\, \,(j\in X)$ form a basis for $\mathcal{E}_G(\lambda)$. Such a set $X$ is called a \textit{star set} for $\lambda$ in $G$ or simply a $\lambda$-star set of $G$. If $X$ is a $\lambda$-star set of $G$ then $\overline{X}=V(G)\setminus X$ is called a $\lambda$-\textit{co-star set} of $G$, while $G-X=G\left[\overline {X}\right]$ is called a \textit{star complement} for $\lambda$ in $G$.

The next result gives some properties of a star set.

\begin{theorem}\label{estrela}\emph{\cite[Proposition 5.1.1]{CvetRowSim2010}}
Given a graph $G$, let $\lambda$ be its eigenvalue with multiplicity $k>0.$ The
following conditions on a vertex subset $X \subset V(G)$ are equivalent:
\begin{enumerate}
\item $X$ is a $\lambda$-star set of $G$;
\item $\mathbb{R}^{n}=\mathcal{E}_{G}(\lambda) \oplus \mathcal{V}$, where
      $\mathcal{V}=\left \langle e_{i}:i\in \overline{X}\right \rangle ;$
\item $\left \vert X\right \vert =k$ and $\lambda$ is not an eigenvalue of $G-X.$
\end{enumerate}
\end{theorem}

It is also worth recalling the following result, known as the Reconstruction Theorem, that states
another characterization of star sets needed in the sequel.

\begin{theorem}\label{rectheorem}\emph{\cite[p.~140]{CvetRowSim2010} }Let $X\subset V(G)$ be
a set of vertices of a graph $G,$ $\overline{X}=V(G)\setminus X$ and assume that $G$
has adjacency matrix $$A_{G}=\begin{bmatrix}
                                         A_{X} & N^{T}\\
                                         N     & C_{\overline{X}}%
                                  \end{bmatrix},$$ where $A_{X}$ and $C_{\overline{X}}$ are the
adjacency matrices of the subgraphs induced by $X$ and $\overline{X},$ respectively. Then $X$ is a
$\lambda$-star set of $G$ if and only if $\lambda$ is not an eigenvalue of $C_{\overline{X}}$ and
\[
A_{X}-\lambda I_{X}=N^{T}\big[C_{\overline{X}}-\lambda I_{\overline{X}}\big]^{-1}N,
\]
where $I_{X}$ and $I_{\overline{X}}$ are respectively the identity matrices of orders $\left \vert X\right \vert$
and $\left \vert \overline{X}\right \vert.$

Furthermore, $\mathcal{E}_G(\lambda)$ is spanned by the vectors%
\begin{equation*}
\begin{bmatrix}%
             \mathbf{y}\\
             -\left(  C_{\overline{X}}-\lambda I_{\overline{X}}\right) ^{-1}N\mathbf{y}
      \end{bmatrix}, \label{rectheorem_spanvecs}%
\end{equation*}
where $\mathbf{y} \in \mathbb{R}^{\left \vert X\right \vert }.$
\end{theorem}

We now prove the following result which will be used in the sequel.

\begin{lemma}\label{star_complement_submatrix_lemma}
Let $G$ be a graph of order $n$, $\lambda \in \sigma(G)$ and $X \subset V(G)$ a $\lambda$-star
set of $G$. The rows of the submatrix
\begin{equation}\label{star_complement_submatrix_1}
\begin{bmatrix}%
              N & C_{\overline{X}}-\lambda I_{\overline{X}}%
       \end{bmatrix}
\end{equation}
span the row space of the matrix $A_G - \lambda I_n$.
\end{lemma}

\begin{proof}
Since $C_{\overline{X}}-\lambda I_{\overline{X}}$ is non-singular, it follows that the $|\overline{X}|$ rows of \eqref{star_complement_submatrix_1}
are linearly independent. Therefore, the result follows since the null space of the matrix $A_G - \lambda I_n$ has dimension $|X|$.
\end{proof}

In the simplex terminology, every square nonsingular submatrix of order $\vert \overline{X}\vert$ of the matrix
\eqref{star_complement_submatrix_1} is called a basic submatrix {and the remaining submatrix is
non-basic. Accordingly, in \eqref{star_complement_submatrix_1} $C_{\overline{X}}-\lambda I_{\overline{X}}$ is basic and
$N$ is non-basic.} Observe that the matrix \eqref{star_complement_submatrix_1}
has $\vert \overline{X}\vert$ rows and $n$ columns. On the other hand, the submatrix $N$ has $\vert X \vert$ columns.
From the next proposition we may conclude that every basic submatrix of the matrix \eqref{star_complement_submatrix_1}
defines a co-star set and vice versa.

\begin{proposition}\emph{\cite{Cardoso_Luz_2014}}\label{equivsssba}
Let $G$ be a graph of order $n$ with at least one edge and $X \subset V(G)$ be a star set for
$\lambda \in \sigma(G)$. Then $X' \subset V(G)$ is a $\lambda$-star set of $G$ if and only if
the submatrix of \eqref{star_complement_submatrix_1} defined by the columns indexed by the vertices
in the $\lambda$-co-star set $\overline{X'}$ is basic, that is, non-singular.
\end{proposition}

Assuming that $G$ has $m$ distinct eigenvalues $\mu_1 \ge\mu_2\ge \dots \ge \mu_m$, where each eigenvalue
$\mu_i$ has multiplicity $k_i$ (and then $\sum_{i=1}^{m}{k_i}=n$), it can be proved
that there is a partition $X_1 \cup X_2\cup \cdots \cup X_m$ of $V(G)$ where each part $X_i$ is a $\mu_i$-star
set (and then has cardinality $k_i$) \cite{rowmain}. This partition is called a \textit{star partition} of $G$.

A vertex subset $D \subset V(G)$ is called a \textit{dominating set} if each vertex in $\overline{D}=V(G) \setminus D$
is adjacent to a vertex of $D$. Following \cite{Slater88}, we say that the dominating set $D$ is a \textit{location dominating set} if
$N_G(u) \cap D \ne N_G(v) \cap D$ whenever $u,v$ are distinct vertices in~$\overline{D}$.
The \textit{domination number} (respectively, \textit{location-domination number}) of $G$ is the least
cardinality of a dominating set (location-dominating set).

\begin{proposition}\emph{\cite{row1994}}\label{domination_th}
Let $X_1 \cup X_2\cup \cdots \cup X_m$ be a star partition of a graph $G$ and suppose that $G$ has no isolated
vertices. Then
\begin{enumerate}
\item for each $i \in \{1, 2, \dots, m\},$ $\overline{X}_i$ is a dominating set for $G$;
\item if  $\mu_i \not \in \{-1,0\}$, then $\overline{X}_i$ is a location-dominating set for $G$.
\end{enumerate}
\end{proposition}

\section{Main and non-main vertices}\label{sec_3}
For a graph $G$, an eigenvalue $\lambda \in \sigma(G)$ and a $\lambda$-star set
$X \subseteq V(G)$, a vertex $v \in X$ is called \textit{$\lambda$-main} (\textit{$\lambda$-non-main}) if $\lambda$ is a
main (non-main) eigenvalue of the subgraph of $G$ induced by $\overline{X} \cup \{v\}.$

Let $B=C_{\overline{X}}-\lambda I_{\overline{X}}$. Multiplying the submatrix \eqref{star_complement_submatrix_1} by $B^{-1}$,
we obtain
\begin{equation*}\label{star_complement_submatrix_2}
\begin{bmatrix}%
              B^{-1}N & I_{\overline{X}}%
       \end{bmatrix}.
\end{equation*}
This matrix contains  the full information about the eigenvectors of $A_G$ afforded by $\lambda$. In fact, the vectors
$$\begin{bmatrix}
              -\textbf{e}_i \\
              B^{-1}N\textbf{e}_i \\
       \end{bmatrix},$$ where $\textbf{e}_i$ is the $i$-th vector of the canonical basis of $\mathbb{R}^{|X|}$,
with $i \in X$, belong to the null space of the matrix \eqref{star_complement_submatrix_1}.
Since this matrix spans the row space of $A_G - \lambda I_n$,
it follows that these vectors also belong to the null space of $A_G - \lambda I_n$. Therefore, the mentioned vectors
are the $\vert X \vert$ linearly independent eigenvectors of $A_G$ associated with the eigenvalue $\lambda$, that is,
they form a basis for  $\mathcal{E}(\lambda)$, and the eigenvalue $\lambda$ is non-main if and only if
$$\textbf{j}^\intercal \begin{bmatrix}
                          -\textbf{e}_i \\
                          B^{-1}N\textbf{e}_i \\
                           \end{bmatrix} = -1 + \textbf{j}_{B}B^{-1}N\textbf{e}_i = 0,$$ holds for all $i \in X$.
                           Accordingly, $\lambda$ is non-main
if and only if
\begin{equation}\label{cond}\textbf{j}^\intercal_{B}B^{-1}N - \textbf{j}_{N}^\intercal = [0, 0, \ldots, 0],\end{equation}
where $\mathbf{j}_{B}$ ($\mathbf{j}_{N}$) is the all-1 vector with a number of entries equal to the cardinality of the co-star set $\overline{X}$ (star set $X$) defined by $B$ ($N$).

According to the definition, a vertex $i \in X$ is $\lambda$-main if
$\textbf{j}^\intercal_{B}B^{-1}N\textbf{e}_i - 1 \ne 0$. Therefore, by considering the simplex tableau associated  with the $\lambda$-star set ($\lambda$-co-star set)
$X$ ($\overline{X}$),
\begin{equation}\label{simplex_tableau}
\begin{tabular}{c|c|c}
                   & $X_N$                                          &\\ \hline
         $X_B$ & $B^{-1}N$                                          &\\ \hline
                   & $\textbf{j}^\intercal_{B}B^{-1}N-\textbf{j}^\intercal_N$ &
\end{tabular},
\end{equation}
where $X_B=\overline{X}$ and $X_N=X$, we deduce that the number of non-zero entries of the last row (usually called the \textit{reduced cost row}) is equal to the number of main vertices of the $\lambda$-star set $X$.

From the previous analysis we obtain the following proposition.

\begin{proposition}\label{proposition_non-main}
For a graph $G$ without isolated vertices and $\lambda \in \sigma(G)$, let $X \subset V(G)$ be a $\lambda$-star set of $G$. The following statements hold:
\begin{enumerate}
\item $\lambda$ is non-main if and only if
      $\textbf{j}^\intercal_{B}B^{-1}N = \textbf{j}^\intercal_N,$ where $B=C_{\overline{X}}-\lambda I_{\overline{X}}$, that is, if and only if all the vertices in $X$ are non-main.
\item Assuming that $\lambda$ is main, the vertex $i \in X$ is main {(non-main)} if and only if the corresponding entry of the reduced cost row of the simplex tableau \eqref{simplex_tableau} is non-zero {(zero)}.
\end{enumerate}
\end{proposition}

	\begin{figure}[ht]
	\centering
	\includegraphics[width=100mm,angle=0]{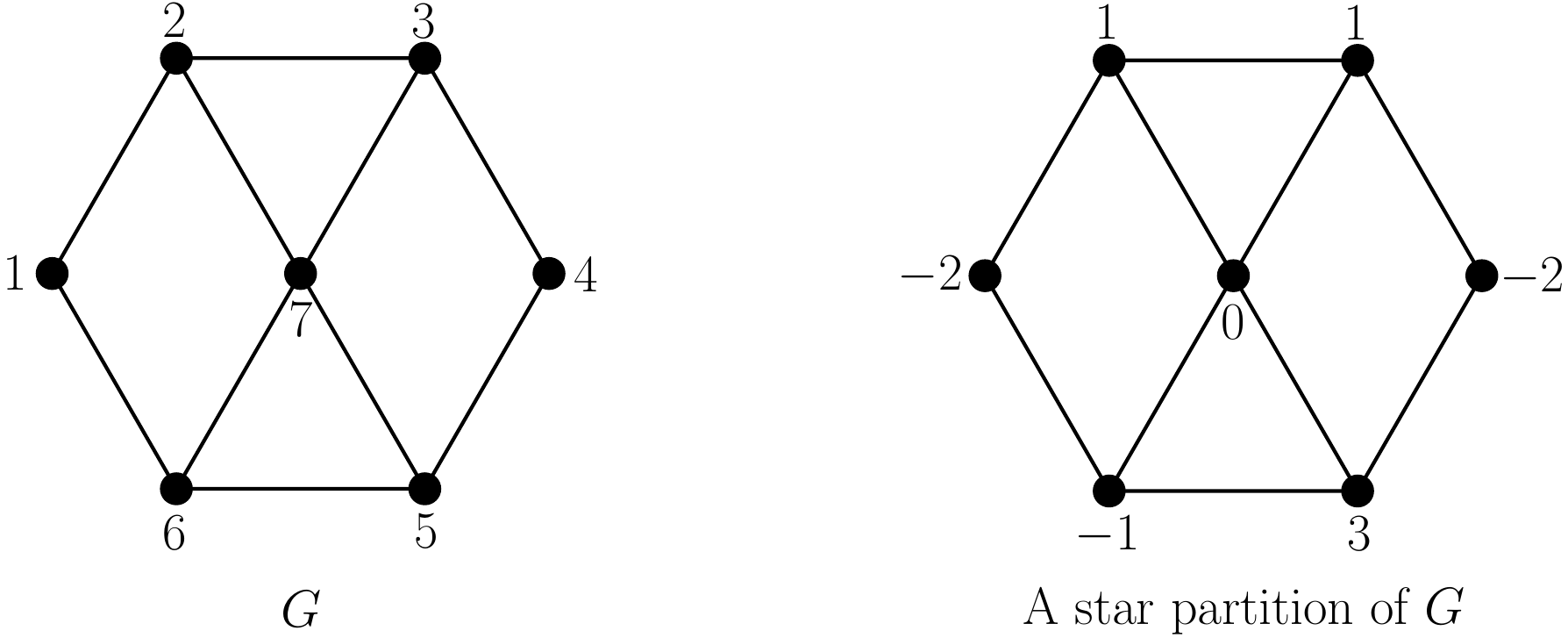}
	\caption{A graph $G$ with $\sigma(G)=\{3, 1^{[2]}, 0, -1,  -2^{[2]}\}$ and a star partition of {$G$}.}\label{figura_1}
\end{figure}

\begin{example}\label{example_1}
Let $G$ be a graph illustrated in Figure~\ref{figura_1}, and let us consider the star sets $X=\{1, 4\}$ and $X'=\{2, 3\}$ of the eigenvalues $-2$ and $1$, respectively.
Then the submatrices $\big[%
                                   N ~~ C_{\overline{X}}-\lambda I_{\overline{X}}%
                            \big]$ and $\big[
                                                                  N ~~ C_{\overline{X'}}-\lambda I_{\overline{X'}}%
                                                           \big]$ are
$$
\begin{bmatrix}
       1 & 0 & 2 & 1 & 0 & 0 & 1\\
       0 & 1 & 1 & 2 & 0 & 0 & 1\\
       0 & 1 & 0 & 0 & 2 & 1 & 1\\
       1 & 0 & 0 & 0 & 1 & 2 & 1\\
       0 & 0 & 1 & 1 & 1 & 1 & 2
\end{bmatrix}\text{ and } \begin{bmatrix}
                                 1 & 0 &-1 & 0 & 0 & 1 & 0\\
                                 0 & 1 & 0 &-1 & 1 & 0 & 0\\
                                 0 & 0 & 0 & 1 &-1 & 1 & 1\\
                                 0 & 0 & 1 & 0 & 1 &-1 & 1\\
                                 1 & 1 & 0 & 0 & 1 & 1 &-1
                           \end{bmatrix},
$$
respectively. In the first matrix, the first two columns correspond to the vertices $1$ and $4$, while the last
$5$ columns correspond to the vertices $2$, $3$, $5$, $6$ and $7$. In the second matrix, the first two columns
correspond to the vertices $2$ and $3$, while the last five columns correspond to the vertices $1$, $4$, $5$, $6$
and $7$. The associated simplex tableaux \eqref{simplex_tableau} are given by
\begin{equation}\label{reduced_tableaus_example1}
\begin{tabular}{r|rr|r}
        $X$   &$x_1$&$x_4$&  \\ \hline
      $x_2$   &  1  & 0   & \\
      $x_3$   &  0  & 1   & \\
      $x_5$   &  0  & 1   & \\
      $x_6$   &  1  & 0   & \\
      $x_7$   & -1  &-1   & \\ \hline
              & $0$ & $0$ &
      \end{tabular} \;\text{ and } \; \begin{tabular}{r|rr|r}
                                          $X'$ &$x_2$&$x_3$&  \\ \hline
                                        $x_1$   &  1  & 0   & \\
                                        $x_4$   &  0  & 1   & \\
                                        $x_5$   &  0  &-1   & \\
                                        $x_6$   & -1  & 0   & \\
                                        $x_7$   &  1  & 1   & \\ \hline
                                                & $0$ & $0$ &
                                  \end{tabular}
\end{equation}
Therefore, by applying Proposition~\ref{proposition_non-main} - item 1, we may conclude that the eigenvalues $-2$ and $1$ are
both non-main {and thus all the vertices in $X$ and $X'$ are non-main}. As another example, considering the eigenvalue $0$, we get that $X''=\{7\}$ is a $0$-star set and then the associated simplex tableau takes the form
\begin{equation}\label{reduced_tableaus_example2}
\begin{tabular}{l|r|r}
        $X''$ &$x_7$&  \\ \hline
      $x_1$   &  1  & \\
      $x_2$   &  0  & \\
      $x_3$   &  0  & \\
      $x_4$   &  1  & \\
      $x_5$   &  0  & \\
      $x_6$   &  0  & \\ \hline
              & $1$ &
      \end{tabular}
\end{equation}
Therefore, the non-zero  cost entry implies that $0$ is a main eigenvalue. Furthermore, it also follows
that $\{1\}, \{4\}$ and $\{7\}$ are the unique $0$-star sets.
\end{example}

The next proposition gives some additional properties of main and non-main vertices.

\begin{proposition}\label{star_set_interchange_prop}
Consider a graph $G$ without isolated vertices, an eigenvalue $\lambda \in \sigma(G)$
and a $\lambda$-star set $X$ of $G$. Let, for $v\in X$, $\textbf{y}_v$ denote the column of
the simplex tableau \eqref{simplex_tableau} associated to $x_v$, that is
\begin{equation}
\textbf{y}_v = \big[C_{\overline{X}}-\lambda I_{\overline{X}}\big]^{-1}\textbf{a}_v, \label{non_basic_column}
\end{equation}
where $\textbf{a}_v$ denotes the column of $N$ in \eqref{star_complement_submatrix_1} corresponding to the vertex $v$.
Then
\begin{enumerate}
\item the set $\{i \in \overline{X}: y_{iv} \ne 0\}$ is non-empty;
\item for every $u \in \{i \in \overline{X}: y_{iv} \ne 0\}$ the following properties hold:
      \begin{enumerate}
      \item the vertex subset $X' = \left(X \setminus \{v\} \right) \cup \{u\}$ is a $\lambda$-star set of $G$;
      \item if $v$ is $\lambda$-main ($\lambda$-non-main), then the vertex $u \in X'$ is $\lambda$-main ($\lambda$-non-main) for $G$.
      \end{enumerate}
\end{enumerate}
\end{proposition}

\begin{proof}
Let us consider simplex tableau \eqref{simplex_tableau} associated to the $\lambda$-star set $X$.
We choose an arbitrary vertex $v \in X$ and consider $\textbf{y}_v $ as in \eqref{non_basic_column}.
\begin{enumerate}
\item Since by Proposition~\ref{domination_th} the vertex set of any $\lambda$-star complement of $G$
      is a dominating set, the vertex $v$ has at least one neighbour in $\overline{X}$, and then $\textbf{a}_v \ne 0$.
      Therefore, there exists at least one entry, say $u$, such that $\textbf{y}_{uv} \ne 0$. Otherwise,
      $\textbf{y}_v=0$ and from \eqref{non_basic_column} it follows that $\textbf{a}_v = 0$, which is
      a contradiction.
\item Now, choose $\textbf{y}_v$ as pivoting column in the simplex tableau \eqref{simplex_tableau}.
      \begin{enumerate}
      \item If the entry $\textbf{y}_{uv}$ is the pivot, then
            $X' = \left(X \setminus \{v\} \right) \cup \{u\}$ is a $\lambda$-star set of $G$.
      \item If the reduced cost in the column associated to $x_v$ is non zero (zero), then the vertex
            $v$ is $\lambda$-main ($\lambda$-non-main) and after pivoting the column associated to $x_u$ remains
            non-zero (zero) and the vertex $u$ becomes $\lambda$-main ($\lambda$-non-main).
      \end{enumerate}
\end{enumerate}
\end{proof}

The next corollary is an immediate consequence of Proposition~\ref{star_set_interchange_prop}.

\begin{corollary}\label{star_complement_cor}
For a graph $G$ without isolated vertices and an eigenvalue $\lambda$, every vertex of $G$
belongs to the vertex set of some $\lambda$-star complement.
\end{corollary}

%\textcolor{red}{Could you explicitly say in the proof of Proposition 9, where you use the assumption that a graph has no isolated vertices?}\textcolor{blue}{ When it is referred that "by Proposition 5 the vertex set of any $\lambda$-star complement of $G$ is a dominating set". Note that in the hypothesis of Proposition 5 it is assumed that $G$ has no isolated vertices. Otherwise, assuming that x is an isolated vertex and then x belongs to some 0-star set X, it is not true that the complement of X is a dominating set}.

The same conclusion can be obtained from \cite[Proposition~7.4.8]{eig} which is proved using a different approach.
From this corollary, and taking into account that every graph admits a star partition, we conclude that there
is more than one star set for every eigenvalue.

A related result is the following.

\begin{proposition}\emph{\cite[Proposition~7.4.8]{eig}}\label{every_star_set}
A vertex $v$ of a graph $G$ lies in every star set corresponding to the eigenvalue $\mu$ if and only if $\mu=0$ and $v$ is an isolated vertex of $G$.
\end{proposition}

Concerning the vertices of some $\lambda$-star set, we prove the following additional result.

\begin{proposition}\label{star_set_vertex_prop.}
Let $G$ be a graph without isolated vertices and $\lambda \in \sigma(G)$. Consider an arbitrary
$\lambda$-star set $X$ and its associated simplex tableau \eqref{simplex_tableau}. A vertex
$v \not \in X$ belongs to some $\lambda$-star set if and only if the row of \eqref{simplex_tableau}
corresponding to $v$ has at least one non-zero entry.
\end{proposition}

\begin{proof}
If there is some non-zero entry $y_{vj}$ in the row of  \eqref{simplex_tableau} corresponding to the
      vertex $v$, from Proposition~\ref{star_set_interchange_prop} - item 2 we get that
      $X' = \left(X \setminus \{j\} \right) \cup \{v\}$ is a $\lambda$-star set of $G$.

Conversely, let us assume that every entry of the row corresponding to the vertex $v$ is zero.
      Considering the submatrix \eqref{star_complement_submatrix_1} and multiplying this submatrix by
      $\big[C_{\overline{X}}-\lambda I_{\overline{X}}\big]^{-1}$, we obtain
      \begin{equation*}\label{matrix_m}
      M = \Big[\big[C_{\overline{X}} - \lambda I_{\overline{X}}\big]^{-1}N \;\;\; I_{\overline{X}}\Big].
      \end{equation*}
      It follows immediately that the row of $M$ assigned to the vertex $v$ has all its entries equal
      to zero, except the diagonal entry of $I_{\overline{X}}$. From the eigenvalue equation, it follows
      that for every $\textbf{u} \in \mathcal{E}_G(\lambda), \; M\textbf{u}=0$ and this equation implies
      $\textbf{u}_v = 0$. Therefore, $v$ does not belong to a {$\lambda$}-star set.
\end{proof}

By virtue of the previous proposition, we deduce that the information available in the simplex tableau \eqref{simplex_tableau} associated to any $\lambda$-star set ($\lambda$-co-star set) of a graph $G$ enables us to detect which vertices have no $\lambda$-star sets. For instance, from the simplex tableau \eqref{reduced_tableaus_example1} of Example~\ref{example_1},
we may conclude that every vertex of the graph of  Figure~\ref{figura_1} belongs to some $1$-star set and
also to some $(-2)$-star set. On the other hand, from the simplex tableau \eqref{reduced_tableaus_example2} we see that there are
only three $0$-star sets: $\{1\}$, $\{4\}$ and $\{7\}$.

\section{Graph parameters related to main and non-main vertices}\label{sec_4}

Given an eigenvalue $\lambda$ of a graph $G$, let  $\main(X)$ denote the subset of $\lambda$-main vertices of the $\lambda$-star set $X$ and $\mathcal{SS}(\lambda,G)$ denote the set of $\lambda$-star sets of $G$. We denote
\begin{align}
\aleph_{\max}(\lambda,G) &=\max \{\vert \main(X) \vert: X \in \mathcal{SS}(\lambda,G)\}, \label{graph_invariant_1}\\
\aleph_{\min}(\lambda,G) &= \min \{\vert \main(X) \vert: X \in \mathcal{SS}(\lambda,G)\}. \label{graph_invariant_2}
\end{align}
Evidently, $\aleph_{\max}(\lambda,G)$ and $\aleph_{\min}(\lambda,G)$ denote the number of $\lambda$-main vertices of a $\lambda$-star set having the maximum number and the minimum number of $\lambda$-main vertices, respectively. One may observe that when, for some $\lambda$-star set $X$, $|\main(X)|=p$ then the number of $\lambda$-non-main vertices in $X$ is equal to $|X|-p$.

Returning to the simplex tableau \eqref{simplex_tableau} and taking into account that the number of main vertices
of the star set $X$ is equal to the number of non-zero entries in the reduced cost row, we can reformulate the
determination of this graph invariant as the determination of the number of non-zero entries in the reduced cost
row of the simplex tableau associated to the basis with maximum number of non-zero entries in the reduced cost row.
{For this purpose, let us define $\delta_1(c)$ as the number of entries equal to $1$ in an arbitrary row vector $c$ and let $(\mathcal{B},\mathcal{N})$ denote the set of partitions of the matrix \eqref{star_complement_submatrix_1} into the pairs of submatrices $(B,N)$ produced by pivoting the associated simplex tableau, where $B$ is basic and $N$ is non-basic. Note that, from the simplex tableau associated to \eqref{star_complement_submatrix_1} by pivoting we may produce all the pairs of basic and non-basic matrices $(B,N)$ of $(\mathcal{B},\mathcal{N})$. Then} the optimization problems \eqref{graph_invariant_1} and \eqref{graph_invariant_2} can be reformulated as follows:
\begin{align*}
\aleph_{\max}(\lambda,G) =&\, k_\lambda-\min\{\delta_1(\mathbf{j}^T_{B}B^{-1}N): (B,N) \in (\mathcal{B},\mathcal{N})\},\\
\aleph_{\min}(\lambda,G) =&\, k_\lambda-\max\{\delta_1(\mathbf{j}^T_{B}B^{-1}N): (B,N) \in (\mathcal{B},\mathcal{N})\},
\end{align*}
where $k_\lambda$ is the multiplicity of the eigenvalue $\lambda$. {Observe that the reduced cost row
of the simplex tableau associated to $(B,N)$, $\mathbf{j}^T_{B}B^{-1}N-\mathbf{j}^T_N$, has the maximum number of non-zero entries when $\mathbf{j}^T_{B}B^{-1}N$ has the minimum number of entries equal to 1}. Therefore, starting from some simplex tableau and using pivot operations we may obtain a sequence of new pairs $(B,N)$ until the above numbers cannot be improved.

In relation to the concepts of $\lambda$-main and $\lambda$-non-main vertices we may define the $\lambda$-main ($\lambda$-non-main) degree of a vertex as follows. Let $G$ be a graph with a main eigenvalue~$\lambda$. The \textit{$\lambda$-main degree} and the \textit{$\lambda$-non-main degree} of a vertex $v \in V(G)$ are
$$\begin{aligned}
d_{(\lambda^+,G)}(v) &= \vert \{S \in \mathcal{SS}(\lambda,G): v \in \main(S)\}\vert,\\
d_{(\lambda^-,G)}(v) &= \vert \{S \in \mathcal{SS}(\lambda,G): v \in S \setminus \main(S)\}\vert,
\end{aligned}$$
respectively. Accordingly, the maximum (minimum) $\lambda$-main degree and the maximum (minimum) $\lambda$-non-main degree of $G$ are
$$ \begin{aligned}
\Delta(\lambda^+,G) (\delta(\lambda^+,G)) = \max(\min)\{d_{(\lambda^+,G)}(v): v \in V(G)\},\\
\Delta(\lambda^-,G) (\delta(\lambda^-,G)) = \max(\min)\{d_{(\lambda^-,G)}(v): v \in V(G)\},
\end{aligned}$$
respectively. As a direct consequence of Proposition~\ref{every_star_set}, if $\lambda \ne 0$ or $G$ has no isolated vertices, then the inequality
$$
d_{(\lambda^+,G)}(v) + d_{(\lambda^-,G)}(v) < |\mathcal{SS}(\lambda,G)|
$$
holds for every vertex $v \in V(G)$.

It is immediate that if $\lambda$ is a main eigenvalue of $G$ and $X$ is a $\lambda$-star
set in which every vertex is main, then $\aleph(\lambda,G)=\vert X \vert$. Furthermore, taking into account
Proposition~\ref{proposition_non-main} - item 2, we get that if  $\lambda \in \sigma(G)$ is main, then
$\aleph_{\min}(\lambda,G) \ge 1$; otherwise, $\aleph_{\max}(\lambda,G)=\aleph_{\min}(\lambda,G)=0$.

The foregoing invariants can be used as tools to check if two graphs are not isomorphic. Namely, the following proposition states several necessary conditions for main eigenvalues of  isomorphic graphs. (We recall that two graphs $G$ and $H$ are isomorphic if and only if there exists a permutation matrix $P$ such that $PA_GP^\intercal = A_H$.)

\begin{proposition}\label{properties_for_ip}
Let $G$ and $H$ be isomorphic graphs. Then they share the same main
eigenvalues. In addition, for each main eigenvalue $\lambda$ the following properties hold.
\begin{enumerate}
\item[(a)] $\vert \mathcal{SS}(\lambda,G)| = |\mathcal{SS}(\lambda,H) \vert$;
\item[(b)] $\aleph_{\max}(\lambda,G) = \aleph_{\max}(\lambda,H) ~~ \text{ and } ~~ \aleph_{\min}(\lambda,G) = \aleph_{\min}(\lambda,H)$;
\item[(c)] $\vert \{X \in \mathcal{SS}(\lambda,G): |\main(X)| = p\} \vert =
                            \vert \{Y \in \mathcal{SS}(\lambda,H): |\main(Y) \vert = p\} \vert,$
                            for $\aleph_{\min}(\lambda,G) \le p \le \aleph_{\max}(\lambda,G)$;
\item[(d)] $\Delta(\lambda^+,G) (\delta(\lambda^+,G))=\Delta(\lambda^+,H) (\delta(\lambda^+,H))$;
\item[(e)] $\Delta(\lambda^-,G) (\delta(\lambda^-,G))=\Delta(\lambda^-,H) (\delta(\lambda^-,H))$;
\item[(f)] $\vert \{v \in V(G): d_{(\lambda^+,G)}(v)=q\} \vert = \vert \{v \in V(H): d_{(\lambda^+,H)}(v)=q\} \vert$,
           for $\delta(\lambda^+,G) \le q \le \Delta(\lambda^+,G)$;
\item[(g)] $\vert \{v \in V(G): d_{(\lambda^-,G)}(v)=q\} \vert = \vert \{v \in V(H): d_{(\lambda^-,H)}(v)=q\} \vert$, for $\delta(\lambda^-,G) \le q \le \Delta(\lambda^-,G)$;
\item[(h)] If $A$ and $B$ are the vertex subsets of $G$ and $H$, respectively, with the same $\lambda$-main ($\lambda$-non-main) degree, then they share the same combinatorial properties as the list of vertex degrees and isomorphic induced subgraphs.
\end{enumerate}
\end{proposition}

\begin{proof}Since none of the considered graph parameters or combinatorial substructures, like vertex degrees and induced subgraphs, changes when the vertices of a graph $G$ are permuted, that is, when its adjacency matrix $A_G$ becomes $P A_GP^\intercal$, where $P$ is a permutation matrix, all the properties follow immediately.
\end{proof}

As it is well-known, the largest eigenvalue of a connected graph is main and simple. The application of Proposition~\ref{properties_for_ip}(a)-(g) to pairs of connected cospectral graphs $G$ and $H$ of order $n$, with the largest eigenvalue in the role of the main eigenvalue $\lambda$ is inconclusive. Indeed, since $\lambda$ decreases when any vertex of $G$ ($H$) is deleted, we have that every vertex of $G$ ($H$) is $\lambda$-main and form a $\lambda$-star set, that is, $\mathcal{SS}(\lambda,G)=\{\{v\}: v \in V(G)\}$ ($\mathcal{SS}(\lambda,H)=\{\{v\}: v \in V(H)\}$). Therefore,
\begin{itemize}
\item[(1)] $\vert \mathcal{SS}(\lambda,G)| = |\mathcal{SS}(\lambda,H) \vert = n$;
\item[(2)] $\aleph_{\max}(\lambda,G) = \aleph_{\max}(\lambda,H) = \aleph_{\min}(\lambda,G) = \aleph_{\min}(\lambda,H) = 1$;
\item[(3)] $\vert \{X \in \mathcal{SS}(\lambda,G): |\main(X)| = 1\} \vert =
                            \vert \{Y \in \mathcal{SS}(\lambda,H): |\main(Y) \vert = 1\} \vert = n$;
\item[(4)] $\Delta(\lambda^+,G) (\delta(\lambda^+,G))=\Delta(\lambda^+,H) (\delta(\lambda^+,H)) = 1$;
\item[(5)] $\Delta(\lambda^-,G) (\delta(\lambda^-,G))=\Delta(\lambda^-,H) (\delta(\lambda^-,H)) = 0$;
\item[(6)] $\vert \{v \in V(G): d_{(\lambda^+,G)}(v)=1\} \vert = \vert \{v \in V(H): d_{(\lambda^+,H)}(v)=1\} \vert = n$;
\item[(7)] $\vert \{v \in V(G): d_{(\lambda^-,G)}(v)=0\} \vert = \vert \{v \in V(H): d_{(\lambda^-,H)}(v)=0\} \vert = n$;
\item[(8)] The vertex subsets $A$ and $B$ of $G$ and $H$, respectively, with the same $\lambda$-main ($\lambda$-non-main) degree are $A=V(G)$ and $B=V(H)$.
\end{itemize}

Consequently, Proposition~\ref{properties_for_ip}(a)-(g) cannot decide whether two cospectral regular graphs are not isomorphic since the necessary conditions are satisfied for every such a pair. The condition (h) may fail for some combinatorial structures, however its verification requires the comparison of the graphs as a whole. We note that the smallest cospectral regular graphs have 10 vertices \cite[p.~10]{StaR}.

\begin{remark}
The analysis of the computational complexity of the determination of all invariants in Proposition~\ref{properties_for_ip} can be done considering the following algorithm.
% Algoritmo
\begin{algorithm}[H]%
	\textbf{Requires:} {A graph $G$ of order $n$ without isolated vertices and a main eigenvalue $\lambda$ with multiplicity $k \ge 1$.}

	\textbf{Ensures:} The graph invariants in (a)--(h) of Proposition~\ref{properties_for_ip}.
	
\begin{algorithmic}[1]
\STATE Determine a star set $X \subset V(G)$ for the main eigenvalue $\lambda$;
\STATE Determine the simplex tableaux associated to $X$;\label{simplexT1}
\STATE Determine $\mathcal{SS}(\lambda,G)$ pivoting a sequence of simplex tableaux, starting with the simplex tableaux of Step \ref{simplexT1} as exemplified in the Appendix;\label{simplex_sequence}
\STATE Determine the parameters as exemplified in Tables~\ref{Table1} and~\ref{Table2} of the Appendix;
\STATE Using the data obtained in the previous step, determine the invariants involved in the properties (a)--(g) of Proposition~\ref{properties_for_ip} as well as the vertex subsets referred in (h), the list of vertex degrees of their induced subgraphs and, as much as possible, other combinatorial properties.
\end{algorithmic}%
\caption{Computation of the graph invariants for Proposition~\ref{properties_for_ip}\label{algorithm}}%
\end{algorithm}%

Let us analyse the complexity of each  step of Algorithm~\ref{algorithm}.

\begin{enumerate}
\item Consider the $n \times k$ matrix $U$ whose columns are the $k$ linearly independent eigenvectors associated with $\lambda$ and compute a matrix $U'$ obtained from $U$ after pivoting operations and  multiplications by scalars until the identity matrix $I_k$ appears as a submatrix of $U'$. Then it is immediate that the vertex
    subset associated to the indices defining $I_k$ is a star set $X$. Since the determination of the mentioned eigenvectors is polynomial, the determination of $X$ is also polynomial.
\item Assuming that the matrix $\begin{bmatrix}%
                                  N & C_{\overline{X}}-\lambda I_{\overline{X}}%
                                \end{bmatrix}$ is the submatrix of $A_G - \lambda I$ obtained after deleting the
      rows with indices associated to the vertices in $X$, the inverse of the basic matrix
      $B=C_{\overline{X}}-\lambda I_{\overline{X}}$ should be determined. The computations of $B^{-1}$ and $B^{-1}N$ are both polynomial. Therefore, the complexity of this step is polynomial.
\item This is the critical step regarding the complexity of the entire algorithm. Indeed, in the worst case, we can have ${n \choose k}$ star sets. However, in practice,  we can deal with a main eigenvalue with multiplicity $1$ (distinct from the largest eigenvalue) or $2$. If $k=1$ then $N$ is an $n \times 1$ matrix
    and the determination of the star sets is immediate from the non-zero entries of the vector $B^{-1}N$. Observe that each star set is a singleton and for every vertex the $\lambda$-non-main degree is zero and the $\lambda$-main degree is $0$ or $1$ (see \eqref{reduced_tableaus_example2} in Example~\ref{example_1}). If $k=2$, then the upper bound on the number of simplex iterations is  the triangular number $t_{n-1}=n(n-1)/2$. In our examples given in the Appendix the number of star sets is much less than this upper bound.
\item This is a easy step with a polynomial complexity.
\item This step is also polynomial, as it deals with the determination of the lists of vertex degrees and several additional combinatorial parameters whose determination has a polynomial complexity.
\end{enumerate}
Therefore, the overall complexity of Algorithm~\ref{algorithm} has the same order as the complexity of Step~\ref{simplex_sequence} and the complexity of this step is polynomial when $k \in \{1,2\}$. Furthermore, this step can be formulated as a combinatorial optimization problem for which new algorithms can be developed along with a deep study of their complexity. Note that, despite the excellent practical performance of the simplex method (the average number of pivot steps is linear \cite{Haimovich} (see also \cite{Adler_et_al})), the search for a strongly polynomial time simplex algorithm remains as one of the most challenging open problems in Optimization and Discrete Geometry \cite{Smale2000}. The number of iterations of a strongly polynomial time algorithm is bounded by a polynomial in the problem dimensions (and not in the size of the input data) \cite{Schrijver}.
\end{remark}

%\bigskip

Now we demonstrate the use of Proposition~\ref{properties_for_ip}.

\begin{example}\label{ex_spectral_tools}
Let $G$ and $H$ be the pair of cospectral graphs depicted in Figure~\ref{figura_5}. These graphs
appear in \cite[Figure~4.3]{CvetRowSim2010}.

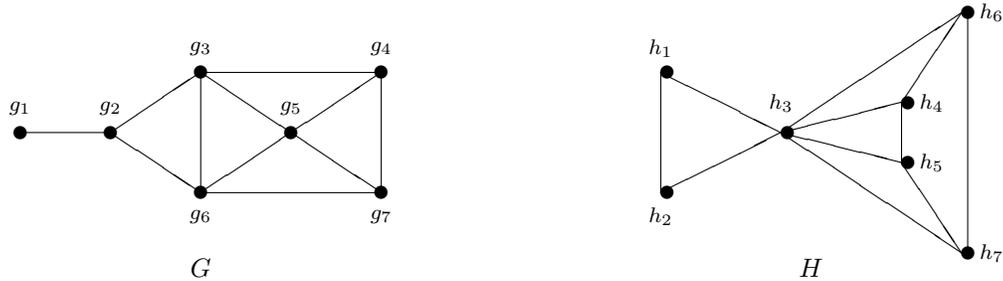
\begin{figure}[h]
\begin{center}
\unitlength=0.4 mm
\begin{picture}(100,80)(-260,70)
%
% vertices do grafo G
\put(-310,130){\circle*{4}} % g3
\put(-250,130){\circle*{4}} % g4
\put(-250,90){\circle*{4}}  % g7
\put(-310,90){\circle*{4}}  % g6
\put(-280,110){\circle*{4}} % g5
\put(-340,110){\circle*{4}} % g2
\put(-370,110){\circle*{4}} % g1
%
% labels do grafo G
\put(-310,138){\makebox(0,0){\footnotesize $g_3$}}
\put(-250,138){\makebox(0,0){\footnotesize $g_4$}}
\put(-250,82){\makebox(0,0){\footnotesize $g_7$}}
\put(-310,82){\makebox(0,0){\footnotesize $g_6$}}
\put(-280,118){\makebox(0,0){\footnotesize $g_5$}}
\put(-340,118){\makebox(0,0){\footnotesize $g_2$}}
\put(-370,118){\makebox(0,0){\footnotesize $g_1$}}
%
%Label of G
\put(-310,65){\makebox(0,0){$G$}}
% edges do grafo G
\put(-310,130){\line(3,-2){60}} %(g_3,g_7)
\put(-310,90){\line(3,2){60}}%(g_4,g_6}
\put(-310,130){\line(1,0){60}} %(g_3,g_4)
\put(-310,90){\line(1,0){60}}%(g_6,g_7}
%\put(-280,110){\line(1,0){30}} %(g_5,g_6)
\put(-310,90){\line(0,1){40}} %(g_6,g_3)
\put(-250,90){\line(0,1){40}} %(g_7,g_4)
\put(-340,110){\line(3,2){30}} %(g_2,g_3)
\put(-340,110){\line(3,-2){30}} %(g_2,g_6)
\put(-370,110){\line(1,0){30}} %(g_1,g_2)
%
% vertices do grafo H
\put(-155,130){\circle*{4}} % h_1
\put(-155,90){\circle*{4}}  % h_2
\put(-115,110){\circle*{4}}% h_3
\put(-75,120){\circle*{4}} % h_4
\put(-75,100){\circle*{4}} % h_5
\put(-55,150){\circle*{4}} % h_6
\put(-55,70){\circle*{4}} % h_7
-
%\put(-40,110){\circle*{4}} % h_2
%\put(-150,110){\circle*{4}} % h_5
%\put(-10,110){\circle*{4}} % h_6
%\put(-10,110){\circle*{4}} % h_5
%
% labels do grafo H
\put(-160,138){\makebox(0,0){\footnotesize $h_1$}}
\put(-160,82){\makebox(0,0){\footnotesize $h_2$}}
\put(-120,120){\makebox(0,0){\footnotesize $h_3$}}
\put(-70,120){\makebox(0,0){\footnotesize $h_4$}}
\put(-70,100){\makebox(0,0){\footnotesize $h_5$}}
\put(-50,150){\makebox(0,0){\footnotesize $h_6$}}
\put(-50,70){\makebox(0,0){\footnotesize $h_7$}}

%
% edges do grafo G
\put(-160,130){\line(2,-1){40}} %(h_1,h_3)
\put(-160,90){\line(2,1){40}}%(h_2,h_3}
\put(-120,110){\line(4,1){40}} %(h_3,h_4)
\put(-120,110){\line(4,-1){40}}%(h_3,h_5}
\put(-160,130){\line(0,-1){40}} %(h_1,h_2)
\put(-80,120){\line(0,-1){20}} %(h_4,h_5)
\put(-120,110){\line(3,2){60}} %(h_3,h_6)
\put(-120,110){\line(3,-2){60}} %(h_3,h_7}
\put(-80,120){\line(2,3){20}} %(h_4,h_6)
\put(-80,100){\line(2,-3){20}} %(h_5,h_7}
\put(-58,150){\line(0,-1){80}} %(h_6,h_7}
%
%Label of H
\put(-110,65){\makebox(0,0){$H$}}
\end{picture}
\medskip
\caption{A pair of cospectral graphs with the common characteristic
         polynomial $p(x)=-16 x^2 - 16 x^3 + 10 x^4 + 11 x^5 - x^7$.} \label{figura_5}
\end{center}
\end{figure}
Taking into account that $0$ is an eigenvalue of $G$ and $H$ with multilpicity $2$, let us consider the $0$-star
sets $X=\{g_6,g_7\}$ and $Y=\{h_6,h_7\}$ of $G$ and $H$, respectively. Then we have
$$
C_{\overline{X}}=\bordermatrix{ & g_1 & g_2 & g_3 & g_4 & g_5 \cr
                            g_1 &  0  &  1  &  0  &  0  &  0  \cr
                            g_2 &  1  &  0  &  1  &  0  &  0  \cr
                            g_3 &  0  &  1  &  0  &  1  &  1  \cr
                            g_4 &  0  &  0  &  1  &  0  &  1  \cr
                            g_5 &  0  &  0  &  1  &  1  &  1  \cr}
~\text{ and }~  C_{\overline{Y}}=\bordermatrix{ & h_1 & h_2 & h_3 & h_4 & h_5 \cr
                                          h_1 &  0  &  1  &  1  &  0  &  0  \cr
                                          h_2 &  1  &  0  &  1  &  0  &  0  \cr
                                          h_3 &  1  &  1  &  0  &  1  &  1  \cr
                                          h_4 &  0  &  0  &  1  &  0  &  1  \cr
                                          h_5 &  0  &  0  &  1  &  1  &  0  \cr}.
$$
The simplex tableaux
\begin{tabular}{c|c|c}
               & $X$                                                                                &\\ \hline
$\overline{X}$ & $C_{\overline{X}}^{-1}N$                                                           &\\ \hline
               & $\textbf{j}^\intercal_{\overline{X}}C_{\overline{X}}^{-1}N-\textbf{j}^\intercal_X$ &
\end{tabular}
and
\begin{tabular}{c|c|c}
               & $Y$                                                                                &\\ \hline
$\overline{Y}$ & $C_{\overline{Y}}^{-1}N$                                                           &\\ \hline
               & $\textbf{j}^\intercal_{\overline{Y}}C_{\overline{Y}}^{-1}N-\textbf{j}^\intercal_Y$ &
\end{tabular}
are
\begin{equation}\label{simplex_tableaux}
\begin{tabular}{r|rr|r}
              &$g_6$&$g_7$&  \\ \hline
      $g_1$   &  1  & -1  & \\
      $g_2$   &  0  &  0  & \\
      $g_3$   &  0  &  1  & \\
      $g_4$   &  1  &  0  & \\
      $g_5$   &  0  &  0  & \\ \hline
              &  1  & -1  &
\end{tabular} \;\text{ and } \; \begin{tabular}{r|rr|r}
                                                &$h_6$&$h_7$&  \\ \hline
                                        $h_1$   &  0  & 0   & \\
                                        $h_2$   &  0  & 0   & \\
                                        $h_3$   &  0  & 0   & \\
                                        $h_4$   &  0  & 1   & \\
                                        $h_5$   &  1  & 0   & \\ \hline
                                                & $0$ & $0$ &
                                  \end{tabular}.
\end{equation}
Therefore, from Proposition~\ref{proposition_non-main} - item 2, we get that $0$ is a main eigenvalue for~$G$
and non-main for $H$ and then, by Proposition~\ref{properties_for_ip}, they are not isomorphic.
\end{example}

From the simplex tableaux \eqref{simplex_tableaux} one may also observe that there are more $0$-star sets in $G$ than in $H$. Note that the simplex tableau associated to the $0$-star set of $H$, $\{h_6,h_7\}$, has only two entries that can be chosen to be the pivoting ones, and thus $H$ has just four $0$-star sets: $\{h_6, h_7\}$, $\{h_4, h_6\}$, $\{h_5, h_7\}$ and $\{h_4, h_5\}$. On the other hand, by pivoting the simplex tableau associated to the $0$-star set of $G$, $\{g_6,g_7\}$, we may produce $8$ distinct $0$-star sets.

From Example~\ref{ex_spectral_tools} we deduce the following remark.

{\begin{remark}\label{main_and_non-main_vertices}
Despite $g_7$ is a main vertex for the $0$-star star $X_1=\{g_6,g_7\}$, we may conclude that $g_7$ is
non-main for the $0$-star set $X_2=\{g_1, g_7\}$.
Indeed, by pivoting the above simplex tableau with the non-zero entry of
the $g_1$-row and $g_6$-column in the role of the pivoting element, we arrive at a tableau in which the entry of the reduced cost row associated to $g_7$ is equal to $0$.
Therefore, a vertex can be non-main for some $\lambda$-star set of a main eigenvalue $\lambda$ and main for
some other $\lambda$-star set.
\end{remark}}

In Example~\ref{ex_spectral_tools} we just deal with two of the three cospectral graphs depicted in
\cite[Figure~4.3]{CvetRowSim2010}. The remaining graph is denoted by $F$ and illustrated in Figure~\ref{fig:F}.
On the basis of computations listed in the Appendix we get the following remark related to $G$ and $F$.

	\begin{figure}
	\centering
	\includegraphics[width=60mm,angle=0]{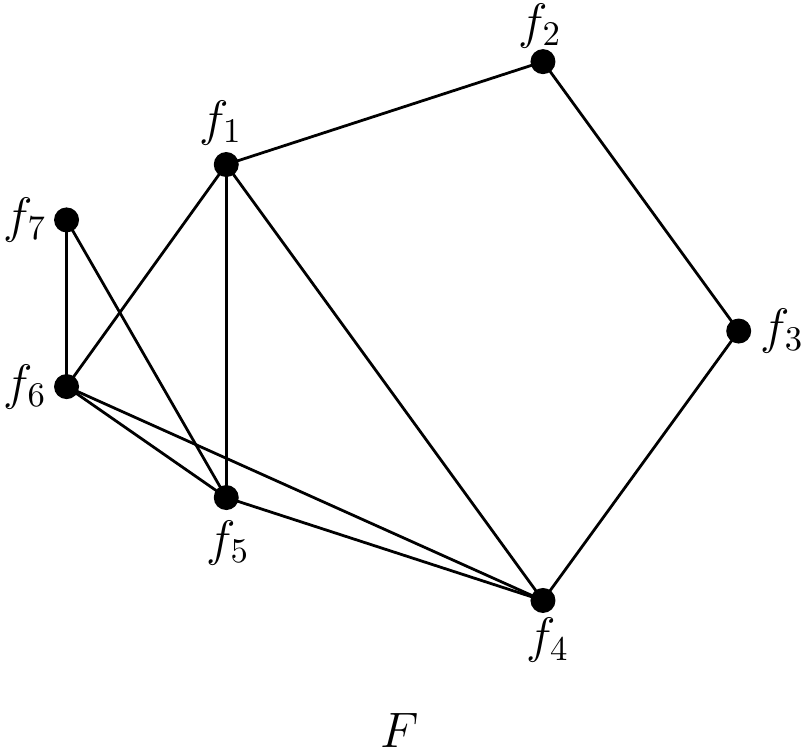}
	\caption{The graph $F$ cospectral with $G$ and $H$ of Example~\ref{ex_spectral_tools}.}\label{fig:F}
\end{figure}

\begin{remark}\label{graphs G_and_F}
It is easy to conclude that $0$ is a main eigenvalue of $G$ and $F$ and the conditions (a)--(g) of Proposition~\ref{properties_for_ip}  hold for this eigenvalue and this pair of graphs. The entire computation of the corresponding parameters, including the list of vertex degrees and induced graphs of subsets of vertices with the same $0$-main ($0$-non-main) degree, is given in the Appendix since it is technical. From the same computation we see that the condition (h) fails to hold, since for instance the vertices $g_1$ of $G$ and $f_7$ of $F$ are the unique vertices in these graphs with $0$-main degree equal to $4$, but on the other hand these vertices differ in degree, which leads to the conclusion that $G$ and $F$ are not isomorphic.
\end{remark}

\section{On maximum value of $\aleph_{\max}(\lambda,G)$}\label{sec:max}

In this section we consider the question of whether  $\aleph_{\max}(\lambda,G)$ is equal to $|X|$, where $X$ is a $\lambda$-star set of $G$.

We know from \cite{row3} that if $H$ is a strongly regular graph with spectrum $\{\nu, \mu^{[k_\mu]}, \lambda^{[k_\lambda]}\}$, where $\nu>\mu>\lambda$, then the cone $K_1\nabla H$ over $H$ has exactly three distinct eigenvalues if and only if $\lambda(\nu-\lambda)=-n$. In this situation, $K_1\nabla H$ has the spectrm $\{\rho, \mu^{[k_\mu]}, \lambda^{[k_\lambda+1]}\}$ and its main eigenvalues are $\rho$ and $\lambda$. The latter can be seen by the fact that $[0, \mathbf{y}^\intercal]^\intercal$ is an eigenvector afforded by $\mu$ in $K_1\nabla H$ if and only if $\mathbf{y}$ is an eigenvector afforded by the same eigenvalue of $H$. Thus, since $\mu$ is non-main in $H$ (as $H$ is regular, it has exactly one main eigenvalue, $\nu$), it is non-main in $K_1\nabla H$ as well, and we conclude that the remaining two eigenvalues must be main (since $K_1\nabla H$ in non-regular, it has more than one main eigenvalue). For an alternative proof the reader is referred to \cite{ACS}. We record this as the following result.

\begin{proposition}
Under the introduced notations, if $H$ is a strongly regular graph with $\lambda(\nu-\lambda)=-n$, then there is a $\lambda$-star set $X$ of $K_1\nabla H$ such that, for every $v\in X$, $\lambda$ is a main eigenvalue of the subgraph induced by $\overline{X}\cup \{v\}$, that is, $\aleph_{\max}(\lambda,K_1\nabla H)=|X|$.
\end{proposition}

\begin{figure}
	\centering
	\includegraphics[width=50mm,angle=0]{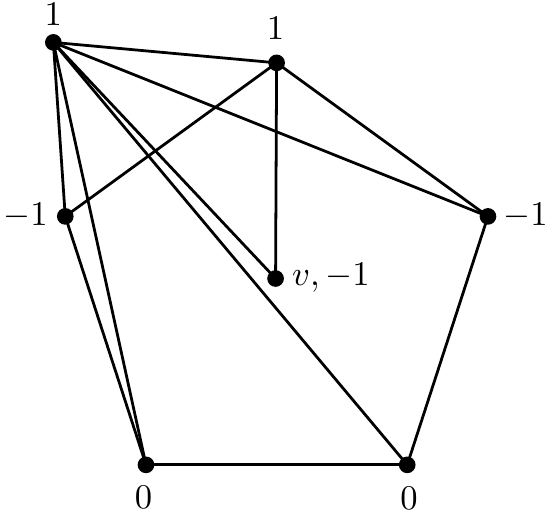}
	\caption{Eigenvector entries for the eigenvalue $\lambda=-2$ in a graph induced by $\overline{X}\cup \{v\}$ of the cone over the Petersen graph.}\label{fig:pet}
\end{figure}

The cone over the Petersen graph serves as an example for the previous proposition. Indeed, the Petersen graph satisfies the equality of the proposition (with $(n, \nu, \lambda)=(10, 3, -2)$), and so $-2$ is a main eigenvalue of the cone. It remains to show that $-2$ is a main eigenvalue for every $G[\overline{X}\cup \{v\}]$, with  a fixed choice of $\overline{X}$ and an arbitrary $v \in X$. The cone over the 5-vertex cycle has no $-2$ as an eigenvalue, so we can take it for the star complement, i.e., the subgraph induced by $\overline{X}$. Then, the subgraphs induced by $\overline{X}\cup\{v\}$, for $v\in X$, are mutually isomorphic as in each of them $v$ is adjacent to exactly two vertices such that exactly one of them belongs to the aforementioned cycle.  So, it is sufficient to consider just one of isomorphic graphs. The eigenvector afforded by $-2$ can be taken to be as in Figure~\ref{fig:pet}, so $-2$ is main, and we are done.

We consider two particular families of graphs in the role of the star complement for an arbitrary eigenvalue $\lambda$: the totally disconnected graphs $tK_1$ and the complete graphs~$K_t$.

\begin{proposition}
	If $tK_1$ is a star complement for an eigenvalue $\lambda~ (\lambda\neq -1)$ in a graph $G$ with $n$ vertices, then $\lambda$ is main in $G$ and $\aleph_{\max}(\lambda,G)=n-t$.
\end{proposition}

\begin{proof}
	We first prove that $\lambda$ is main in $G$. For this purpose, we need to check the equality of \eqref{cond}. In our case, $B^{-1}=-\frac{1}{\lambda}I_t$. Observe that $\lambda\neq 0$, since $\lambda$ does not belong to the spectrum of the star complement. If $\mathbf{n}_v$ is the column of $N$ that corresponds to the vertex $v$ of $X$, then the equality of \eqref{cond} reads $-\frac{1}{\lambda}\mathbf{j}_t\cdot \mathbf{n}_v-1=0$, where $\cdot$ stands for the standard inner product. Obviously, this equality holds if and only if $v$ has exactly $-\lambda$ neighbours in $tK_1$. Since $\lambda$ is an eigenvalue of  $G[\overline{X}\cup\{v\}]$ and the non-zero eigenvalues of this graph are the positive and the negative square root of the number of neighbours of $v$ in $\overline{X}$, we conclude that the equality holds precisely if $\lambda=-1$, the case eliminated in the formulation of the statement. Therefore, since \eqref{cond} does not hold for $v\in X$, we conclude that $\lambda$ is main in $G$.
	
The fact that $\lambda$ is main in  $G[\overline{X}\cup\{v\}]$, for every $v\in X$, is proved in essentially the same way since the only difference is that, in this case, the equality of \eqref{cond} should be checked for a 1-vertex extension of $tK_1$ instead of the entire graph $G$. From \eqref{graph_invariant_1} we get $\aleph_{\max}(\lambda,G)=|X|=n-t$.
\end{proof}

In other words, the equality $\aleph_{\max}(\lambda,G)=|X|$ is attained whenever $\lambda\neq -1$.

The previous proposition is relevant for a negative $\lambda$. For example, by taking $t=8$, $t=10$ and $t=12$, we arrive at the unique maximal graph with $tK_1$ in the role of the star complement for $\lambda=-2$. The first has the spectrum $\{14, 2^{[7]}, -2^{[14]}\}$, and vertex degrees 7 and 16. The second has the spectrum $\{8, 3^{[4]}, 0^{[5]}, -2^{[10]}\}$, and vertex degrees 4 and 10. The third has the spectrum $\{10, 4^{[5]}, 0^{[6]}, -2^{[15]}\}$, and vertex degrees 5 and 12. The first graph is an example of a non-regular graph with exactly 3 distinct eigenvalues.

We proceed with the next result.

\begin{proposition}
If $K_t~(t\geq 2)$ is a star complement for a main eigenvalue $\lambda~(\lambda\neq 0)$ in a graph $G$ with $n$ vertices, then  $\aleph_{\max}(\lambda,G)=n-t$.
\end{proposition}

\begin{proof} Observe that the statement holds for $\lambda=t$, as in this case $G$ is necessarily $K_{t+1}$, i.e., a 1-vertex extension of $K_t$. Observe also that, under the assumption that $t\geq 2$, we have $\lambda\neq -1$, since $\lambda$ does not belong to the spectrum of the star complement.
	
	We need to prove that $\lambda$ is main in the graph induced by $\overline{X}\cup\{v\}$, for every $v\in X$. Suppose that $v$ is adjacent to exactly $p~(p<t)$ vertices of $K_t$. An eigenvector $\mathbf{y}$ afforded by $\lambda$ in the corresponding 1-vertex extension of $K_t$ has at most three distinct coordinates: $y_v$ (that corresponds to $v$) $y'$ (that corresponds to the neighbours of $v$) and $y''$ (that corresponds to non-neighbours of $v$). The eigenvalue equations for $v$ and one of its neighbours are
	\begin{align}\label{eq1}\lambda y_v&=py',\\ \label{eq2} \lambda y'&=y_v+(p-1)y'+(t-p)y'', \end{align}
respectively. If $\lambda$ is non-main, we also have
\begin{equation}\label{eq3}y_v+py'+(t-p)y''=0.\end{equation}	
	From \eqref{eq1} and \eqref{eq3} we get $y'=\frac{\lambda}{p}y_v$ and $y''=\frac{\lambda+1}{p-t}y_v$. Substituting for $y', y''$ in \eqref{eq2}, we arrive at $y_v(\lambda(\lambda+1))=0$. Since $\lambda\notin\{-1, 0\}$, we have $y_v=0$ but this leads to the conclusion that $\mathbf{y}$ is a zero-vector, which is impossible. Hence $\lambda$  is main in $G[\overline{X}\cup\{v\}]$, and we are done.
%
%Now, if $\lambda$ is non-main, we have
%	$$y_v+py'+(t-p)y''=0,$$
%	which together with the previous equality and the assumption that $\lambda\neq 0$ leads to $(py'+(t-p)y'')(1+\frac{1}{\lambda})=0$. Since $\lambda\neq -1$, we have \begin{equation}\label{lam}py'+(t-p)y'=0.\end{equation}
%	Now, the eigenvalue equation for a non-neighbour of $v$ reads $\lambda y''= (t-p-1)y''+py'$, which together with \eqref{lam} gives $\lambda=-1$, a contradiction.
\end{proof}

It is proved in \cite{Sta} that, apart from $K_1$, exactly two complete graphs may appear as star complements for $1$, and then 1 is necessarily the second largest eigenvalue in their extensions. These graphs are $K_{10}$ and $K_{11}$. Moreover, there are exactly two maximal extensions of the former graph. The first has the spectrum $\{11, 1^{[10]}, -1^{[5]}, -4^{[4]}\}$, and vertex degrees 7 and 13. The second one has the spectrum $\{11.28, 1^{[14]}, -1, -3^{[7]}, -3.28\}$, and vertex degrees 5, 9 and 16. On the basis of \eqref{cond}, we confirm that in both $1$ is a main eigenvalue, so these graphs are examples for the previous proposition.

\section{Open problems}\label{sec:problems}

Here we list some conclusions and open problems we spotted during the research. Consider a graph $G$ and a main eigenvalue $\lambda$ of $G$.

\begin{enumerate}
\item From the Appendix we may conclude that, in general, a vertex $v \in V(G)$ which is $\lambda$-main ($\lambda$-non-main) for every $\lambda$-star set may or may not exist. Are there some conditions that would preserve the existence of such a vertex?

\item Example~\ref{ex_spectral_tools} shows that there are vertices $v \in V(G)$ for which there are no $\lambda$-star sets $X$ such that $v$ is $\lambda$-main ($\lambda$-non-main) for $X$. Under which conditions this would be false?

\item{The graphs of Example~\ref{ex_spectral_tools}, where the determination of the introduced new graph invariants is illustrated, are non-isomorphic because they do not share the same vertex degrees. However, in the context of the isomorphism problem, Proposition~\ref{properties_for_ip} should be used in case of cospectral graphs with the same degree sequence. As we noted upon the proposition, it gives only necessary conditions for main eigenvalues of isomorphic graphs, and as already noted the items (a)-(g) are indecisive in the case of cospectral regular graphs. An intriguing problem that arises is to determine the smallest pair of non-regular cospectral graphs with the same degree sequence for which the proposition is indecisive. In relation to this, we can add that our computer search has not find any such a pair with at most 7 vertices.}

\item What is the maximum value of $\aleph_{\max}(\lambda,G)$  among the connected graphs $G$ of order~$n$? Clearly, it is bounded by $|X|$, and according to \cite[Theorem 5.3.1]{CvetRowSim2010}, $|X|$ cannot exceed $t\choose 2$ where $t~(t\geq 2)$ is the codimension of the eigenspace of $\lambda$. In the previous section we have seen some examples of a comparatively large value of $\aleph_{\max}(\lambda,G)$. In fact, in each of these examples $\aleph_{\max}(\lambda,G)$ attains $|X|$, but $|X|$ does not attain its upper bound. So, determining  a sharp upper bound for $\aleph_{\max}(\lambda,G)$ sounds as an interesting research problem.
\end{enumerate}

\section{Appendix}

In what follows we present the computation of the parameters, vertex degrees and induced subgraphs of $G$ and $F$ referred in Remark~\ref{graphs G_and_F} and  Proposition~\ref{properties_for_ip}.

\subsection{The computations for the graph $G$ depicted in Figure~\ref{figura_5}}
Consider the $0$-star $X_1=\{g_6, g_7\}$,
\begin{eqnarray*}
A_G = \bordermatrix{    & g_1 & g_2 & g_3 & g_4 & g_5 & g_6 & g_7 \cr
                    g_1 &  0  &  1  &  0  &  0  &  0  &  0  &  0  \cr
                    g_2 &  1  &  0  &  1  &  0  &  0  &  1  &  0  \cr
                    g_3 &  0  &  1  &  0  &  1  &  1  &  1  &  0  \cr
                    g_4 &  0  &  0  &  1  &  0  &  1  &  0  &  1  \cr
                    g_5 &  0  &  0  &  1  &  1  &  0  &  1  &  1  \cr
                    g_6 &  0  &  1  &  1  &  0  &  1  &  0  &  1  \cr
                    g_7 &  0  &  0  &  0  &  1  &  1  &  1  &  0  \cr} & \text{ and } &
N_G = \bordermatrix{    & g_6 & g_7 \cr
                        g_1 &  0  &  0  \cr
                        g_2 &  1  &  0  \cr
                        g_3 &  1  &  0  \cr
                        g_4 &  0  &  1  \cr
                        g_5 &  1  &  1  \cr}.
\end{eqnarray*}

We have

\begin{eqnarray*}
C_{\overline{X}_1}^{-1} = \bordermatrix{    & g_1 & g_2 & g_3 & g_4 & g_5 \cr
                                       g_1 &-1/2&  1  & 1/2 &-1/2 &-1/2 \cr
                                       g_2 &  1 &  0  &  0  &  0  &  0  \cr
                                       g_3 & 1/2&  0  &-1/2 & 1/2 & 1/2 \cr
                                       g_4 &-1/2&  0  & 1/2 &-1/2 & 1/2 \cr
                                       g_5 &-1/2&  0  & 1/2 & 1/2 &-1/2 \cr} & \text{ and~ }
C_{\overline{X}_1}^{-1}N_G = \bordermatrix{   & g_6 & g_7 \cr
                                         g_1&  1  & -1  \cr
                                         g_2&  0  &  0  \cr
                                         g_3&  0  &  1  \cr
                                         g_4&  1  &  0  \cr
                                         g_5&  0  &  0  \cr}.
\end{eqnarray*}

Here is a sequence of simplex tableaux obtained by pivoting in the one defined by $X_1$ ($\overline{X}_1$).
Each pivoting element appears marked by a framebox.

\medskip

\begin{tabular}{c|c|c}
               & $X_1$                                                                                     &\\ \hline
$\overline{X_1}$ & $C_{\overline{X}_1}^{-1}N$                                                              &\\ \hline
               & $\textbf{j}^\intercal_{\overline{X}_1}C_{\overline{X}_1}^{-1}N_G-\textbf{j}^\intercal_{X_1}$&
\end{tabular} \qquad = \qquad \begin{tabular}{r|rr|r}
                       &$g_6$&$g_7$&  \\ \hline
                 $g_1$ &  \framebox{1}  & -1  & \\
                 $g_2$ &  0  &  0  & \\
                 $g_3$ &  0  &  1  & \\
                 $g_4$ &  1  &  0  & \\
                 $g_5$ &  0  &  0  & \\ \hline
                       &  1  & -1  &
                \end{tabular} \qquad $\rightarrow$ \qquad \begin{tabular}{r|rr|r}
                                                                 &$g_1$&$g_7$&  \\ \hline
                                                           $g_6$ &  1  & -1  & \\
                                                           $g_2$ &  0  &  0  & \\
                                                           $g_3$ &  0  &  1  & \\
                                                           $g_4$ & \framebox{-1}  &  1  & \\
                                                           $g_5$ &  0  &  0  & \\ \hline
                                                                 & -1  &  0  &
                                                           \end{tabular}\\

$\rightarrow$ \;\;\; \begin{tabular}{r|rr|r}
                       &$g_4$&$g_7$&  \\ \hline
                 $g_6$ &  1  &  0  & \\
                 $g_2$ &  0  &  0  & \\
                 $g_3$ &  0  &  \framebox{1}  & \\
                 $g_1$ & -1  & -1  & \\
                 $g_5$ &  0  &  0  & \\ \hline
                       & -1  & -1  &
                \end{tabular} \qquad $\rightarrow$ \qquad \begin{tabular}{r|rr|r}
                                                                 &$g_4$&$g_3$&  \\ \hline
                                                           $g_6$ &  1  &  0  & \\
                                                           $g_2$ &  0  &  0  & \\
                                                           $g_7$ &  0  & 1  & \\
                                                           $g_1$ & -1  & \framebox{1}  & \\
                                                           $g_5$ &  0  &  0  & \\ \hline
                                                                 & -1  &  1  &
                                                           \end{tabular}
                \qquad $\rightarrow$ \qquad \begin{tabular}{r|rr|r}
                                                   &$g_4$&$g_1$&  \\ \hline
                                             $g_6$ &\framebox{1}  & 0  & \\
                                             $g_2$ &  0  &  0  & \\
                                             $g_7$ &  1  & -1  & \\
                                             $g_3$ & -1  &  1  & \\
                                             $g_5$ &  0  &  0  & \\ \hline
                                                   &  0  & -1  &
                                             \end{tabular}

$\rightarrow$ \;\;\; \begin{tabular}{r|rr|r}
                       &$g_6$&$g_1$&  \\ \hline
                 $g_4$ &  1  &  0  & \\
                 $g_2$ &  0  &  0  & \\
                 $g_7$ & -1  & -1  & \\
                 $g_3$ &  1  &\framebox{1}  & \\
                 $g_5$ &  0  &  0  & \\ \hline
                       &  0  & -1  &
                \end{tabular} \qquad $\rightarrow$ \qquad \begin{tabular}{r|rr|r}
                                                                 &$g_6$&$g_3$&  \\ \hline
                                                           $g_4$ &  1  &  0  & \\
                                                           $g_2$ &  0  &  0  & \\
                                                           $g_7$ &  0  &  1  & \\
                                                           $g_1$ &\framebox{1}  & 1 & \\
                                                           $g_5$ &  0  &  0  & \\ \hline
                                                                 &  1  &  1  &
                                                           \end{tabular}
                \qquad $\rightarrow$ \qquad \begin{tabular}{r|rr|r}
                                                   &$g_1$&$g_3$&  \\ \hline
                                             $g_4$ & -1  & -1  & \\
                                             $g_2$ &  0  &  0  & \\
                                             $g_7$ &  0  &  1  & \\
                                             $g_6$ &  1  &  1  & \\
                                             $g_5$ &  0  &  0  & \\ \hline
                                                   & -1  &  0  &
                                             \end{tabular}.

\medskip

It follows that the $0$-star sets of $G$ are the vertex subsets
$X_1=\{g_6, g_7\}$ (with $\main(X_1)=X_1$), $X_2=\{g_1, g_7\}$ ($\main(X_2)=\{g_1\}$), $X_3=\{g_4, g_7\}$
($\main(X_3)=X_3$),  $X_4=\{g_3, g_4\}$ ($\main(X_4)=X_4$), $X_5=\{g_4, g_1\}$ ($\main(X_5)=\{g_1\}$),
$X_6=\{g_1, g_6\}$ ($\main(X_6)=\{g_1\}$), $X_7=\{g_3, g_6\}$ ($\main(X_7)=X_7$) and $X_8=\{g_1, g_3\}$
($\main(X_8)=\{g_1\}$). Table\ref{Table1} summarizes the elements of each $0$-star set and gives the main and
non-main degrees. Each entry $(g_i,X_j)$ is equal to $\left\{\begin{array}{rl}
                                            1 & \text{if~}~ g_i \in \main(X_j),\\
                                           -1 & \text{if~}~ g_i \notin \main(X_j).
                                          \end{array}\right.$\\
\begin{table}
	\caption{Computation of invariants for the main $0$-star sets of $G$}\label{Table1}
\begin{center}
{\small \begin{tabular}{l|c|c|c|c|c|c|c|c|c|c}
      & $X_1$ & $X_2$ & $X_3$ & $X_4$ & $X_5$ & $X_6$ & $X_7$ & $X_8$ & main degree & non-main degree\\ \hline
$g_1$ &       &   1   &       &       &   1   &   1   &       &   1   &   4         &   0 \\
$g_2$ &       &       &       &       &       &       &       &       &   0         &   0 \\
$g_3$ &       &       &       &   1   &       &       &   1   &  -1   &   2         &   1 \\
$g_4$ &       &       &   1   &   1   &  -1   &       &       &       &   2         &   1 \\
$g_5$ &       &       &       &       &       &       &       &       &   0         &   0 \\
$g_6$ &   1   &       &       &       &       &  -1   &   1   &       &   2         &   1 \\
$g_7$ &   1   &  -1   &   1   &       &       &       &       &       &   2         &   1 \\ \hline
\end{tabular}}
\end{center}
\end{table}

Using the obtained data we obtain the following parameters where the itemization refers to that of Proposition~\ref{properties_for_ip}.

\begin{itemize}
\item[(a)] $|\mathcal{SS}(0,G)|=8$;
\item[(b)] $\aleph_{\max}(0,G)=2$ and $\aleph_{\min}(0,G)=1$;
\item[(c)] $|\{X \in \mathcal{SS}(0,G): |\main(X)|=1\}|=4$ and $|\{X \in \mathcal{SS}(0,G): |\main(X)|=2\}|=4$;
\item[(d)] $\delta(0^+,G) = 0$ and $\Delta(0^+,G) = 4$;
\item[(e)] $\delta(0^-,G) = 0$ and $\Delta(0^-,G) = 1$;
\item[(f)] \begin{eqnarray*}
           \vert \{v \in V(G): d_{(0^+,G)}(v)=0\} \vert &=& 2, \\
           \vert \{v \in V(G): d_{(0^+,G)}(v)=1\} \vert &=& 0, \\
           \vert \{v \in V(G): d_{(0^+,G)}(v)=2\} \vert &=& 4, \\
           \vert \{v \in V(G): d_{(0^+,G)}(v)=3\} \vert &=& 0, \\
           \vert \{v \in V(G): d_{(0^+,G)}(v)=4\} \vert &=& 1;
           \end{eqnarray*}
\item[(g)] \begin{eqnarray*}
            \vert \{v \in V(G): d_{(0^-,G)}(v)=0\} \vert &=& 3,\\
            \vert \{v \in V(G): d_{(0^-,G)}(v)=1\} \vert &=& 4;
           \end{eqnarray*}
\item[(h)] Let $V^+_d$ and $V^-_{d}$ be, respectively, the subsets of vertices with $0$-main degree and $0$-non-main degree equal to $d$.
           \begin{enumerate}
           \item $V^+_0=\{g_2, g_5\}$ is an independent set; $d_G(g_2)=3$ and $d_G(g_5)=4$.
           \item $V^+_2=\{g_3, g_4, g_6, g_7\}$; the induced subgraph $G[V^+_2]$ is isomorphic to the cycle
                 $C_4$; $d_{G}(g_3)=d_{G}(g_6)=4$ and $d_G(g_4)=d_G(g_7)=3$.
           \item $V^+_4=\{g_1\}$; $d_G(g_1)=1$.
           \item $V^-_0=\{g_1, g_2, g_5\}$; the induced subgraph $G[V^-_0]$ is isomorphic to $K_1 \cup K_2$; $d_G(g_1)=1$, $d_G(g_2)=3$ and $d_G(g_5)=4$.
           \item $V^-_1=V^+_2$.
           \end{enumerate}
\end{itemize}

\subsection{The computations for the graph $F$ depicted in Figure~\ref{fig:F}}
Consider the $0$-star set of $Y_1=\{f_4, f_7\}$,
\begin{eqnarray*}
A_F = \bordermatrix{    & f_1 & f_2 & f_3 & f_4 & f_5 & f_6 & f_7 \cr
                    f_1 &  0  &  1  &  0  &  1  &  1  &  1  &  0  \cr
                    f_2 &  1  &  0  &  1  &  0  &  0  &  0  &  0  \cr
                    f_3 &  0  &  1  &  0  &  1  &  0  &  0  &  0  \cr
                    f_4 &  1  &  0  &  1  &  0  &  1  &  1  &  0  \cr
                    f_5 &  1  &  0  &  0  &  1  &  0  &  1  &  1  \cr
                    f_6 &  1  &  0  &  0  &  1  &  1  &  0  &  1  \cr
                    f_7 &  0  &  0  &  0  &  0  &  1  &  1  &  0  \cr} & \text{ and } &
N_F = \bordermatrix{    & f_4 & f_7 \cr
                    f_1 &  1  &  0  \cr
                    f_2 &  0  &  0  \cr
                    f_3 &  1  &  0  \cr
                    f_5 &  1  &  1  \cr
                    f_6 &  1  &  1  \cr}.
\end{eqnarray*}

As before, we get

\begin{eqnarray*}
C_{\overline{Y}_1}^{-1} = \bordermatrix{    & f_1 & f_2 & f_3 & f_5 & f_6 \cr
                                       f_1 &-1/2&  0  & 1/2 & 1/2 & 1/2 \cr
                                       f_2 &  0 &  0  &  1  &  0  &  0  \cr
                                       f_3 & 1/2&  1  &-1/2 &-1/2 &-1/2 \cr
                                       f_5 & 1/2&  0  &-1/2 &-1/2 & 1/2 \cr
                                       f_6 & 1/2&  0  &-1/2 & 1/2 &-1/2 \cr}, & \text{  ~}
C_{\overline{Y}_1}^{-1}N_F = \bordermatrix{   & f_4 & f_7 \cr
                                         f_1&  1  &  1  \cr
                                         f_2&  1  &  0  \cr
                                         f_3& -1  & -1  \cr
                                         f_5&  0  &  0  \cr
                                         f_6&  0  &  0  \cr}
\end{eqnarray*}
and

\medskip

\begin{tabular}{c|c|c}
               & $Y_1$                                                                                     &\\ \hline
$\overline{Y_1}$ & $C_{\overline{Y}_1}^{-1}N$                                                              &\\ \hline
               & $\textbf{j}^\intercal_{\overline{Y}_1}C_{\overline{Y}_1}^{-1}N_F-\textbf{j}^\intercal_{Y_1}$&
\end{tabular} \qquad = \qquad \begin{tabular}{r|rr|r}
                       &$f_4$&$f_7$&  \\ \hline
                 $f_1$ &  1  &  1  & \\
                 $f_2$ &  1  &  0  & \\
                 $f_3$ & \framebox{-1}  & -1  & \\
                 $f_5$ &  0  &  0  & \\
                 $f_6$ &  0  &  0  & \\ \hline
                       &  0  & -1 &
                \end{tabular} \qquad $\rightarrow$ \qquad \begin{tabular}{r|rr|r}
                                                                 &$f_3$&$f_7$&  \\ \hline
                                                           $f_1$ &  1  &  0  & \\
                                                           $f_2$ &\framebox{1}  & -1  & \\
                                                           $f_4$ & -1  &  1  & \\
                                                           $f_5$ &  0  &  0  & \\
                                                           $f_6$ &  0  &  0  & \\ \hline
                                                                 &  0  & -1  &
                                                           \end{tabular}\\

$\rightarrow$ \;\;\; \begin{tabular}{r|rr|r}
                       &$f_2$&$f_7$&  \\ \hline
                 $f_1$ &\framebox{-1}  &  1  & \\
                 $f_3$ &  1  & -1  & \\
                 $f_4$ &  1  &  0 & \\
                 $f_5$ &  0  &  0  & \\
                 $f_6$ &  0  &  0  & \\ \hline
                       &  0  & -1  &
                \end{tabular} \qquad $\rightarrow$ \qquad \begin{tabular}{r|rr|r}
                                                                 &$f_1$&$f_7$&  \\ \hline
                                                           $f_2$ & -1  & -1  & \\
                                                           $f_3$ &  1  &  0  & \\
                                                           $f_4$ &  1  &\framebox{1}  & \\
                                                           $f_5$ &  0  &  0  & \\
                                                           $f_6$ &  0  &  0  & \\ \hline
                                                                 &  0  & -1  &
                                                           \end{tabular}
                \qquad $\rightarrow$ \qquad \begin{tabular}{r|rr|r}
                                                   &$f_1$&$f_4$&  \\ \hline
                                             $f_2$ &  0  &  1  & \\
                                             $f_3$ &\framebox{1}  &  0  & \\
                                             $f_7$ &  1  &  1  & \\
                                             $f_5$ &  0  &  0  & \\
                                             $f_6$ &  0  &  0  & \\ \hline
                                                   &  1  &  1  &
                                             \end{tabular}

$\rightarrow$ \;\;\; \begin{tabular}{r|rr|r}
                       &$f_3$&$f_4$&  \\ \hline
                 $f_2$ &  0  &\framebox{1}  & \\
                 $f_1$ &  1  &  0  & \\
                 $f_7$ & -1  &  1  & \\
                 $f_5$ &  0  &  0  & \\
                 $f_6$ &  0  &  0  & \\ \hline
                       & -1  &  1  &
                \end{tabular} \qquad $\rightarrow$ \qquad \begin{tabular}{r|rr|r}
                                                                 &$f_3$&$f_2$&  \\ \hline
                                                           $f_4$ &  0  &  1  & \\
                                                           $f_1$ &\framebox{1}  &  0  & \\
                                                           $f_7$ & -1  & -1  & \\
                                                           $f_5$ &  0  &  0  & \\
                                                           $f_6$ &  0  &  0  & \\ \hline
                                                                 & -1  & -1  &
                                                           \end{tabular}
                \qquad $\rightarrow$ \qquad \begin{tabular}{r|rr|r}
                                                   &$f_1$&$f_2$&  \\ \hline
                                             $f_4$ &  0  &  1  & \\
                                             $f_3$ &  1  &  0  & \\
                                             $f_7$ &  1  & -1  & \\
                                             $f_5$ &  0  &  0  & \\
                                             $f_6$ &  0  &  0  & \\ \hline
                                                   &  1  & -1  &
                                             \end{tabular}.
                                         \medskip

It follows that the $0$-star sets of $F$ are the vertex subsets
$Y_1=\{f_4, f_7\}$ (with $\main(Y_1)=\{f_7\}$), $Y_2=\{f_3, f_7\}$ ($\main(Y_2)=\{f_7\}$), $Y_3=\{f_2, f_7\}$
($\main(Y_3)=\{f_7\}$),  $Y_4=\{f_1, f_7\}$ ($\main(Y_4)=\{f_7\}$), $Y_5=\{f_1, f_4\}$ ($\main(Y_5)=Y_5$),
$Y_6=\{f_3, f_4\}$ ($\main(Y_6)=Y_6$), $Y_7=\{f_2, f_3\}$ ($\main(Y_7)=Y_7$) and $Y_8=\{f_1, f_2\}$ ($\main(Y_8)=Y_8$). As before, Table~\ref{Table2} summarizes the elements of each $0$-star set and gives the main and non-main degrees.

\begin{table}
	\caption{Computation of invariants for the main $0$-star sets of $F$}\label{Table2}
\begin{center}
{\small \begin{tabular}{l|c|c|c|c|c|c|c|c|c|c}
      & $Y_1$ & $Y_2$ & $Y_3$ & $Y_4$ & $Y_5$ & $Y_6$ & $Y_7$ & $Y_8$ & main degree & non-main degree\\ \hline
$f_1$ &       &       &       &  -1   &   1   &       &       &   1   &   2         &   1 \\
$f_2$ &       &       &  -1   &       &       &       &   1   &   1   &   2         &   1 \\
$f_3$ &       &  -1   &       &       &       &   1   &   1   &       &   2         &   1 \\
$f_4$ &  -1   &       &       &       &   1   &   1   &       &       &   2         &   1 \\
$f_5$ &       &       &       &       &       &       &       &       &   0         &   0 \\
$f_6$ &       &       &       &       &       &       &       &       &   0         &   0 \\
$f_7$ &   1   &   1   &   1   &   1   &       &       &       &       &   4         &   0 \\ \hline
\end{tabular}}
\end{center}
\end{table}

Further, we obtain the following.

\begin{itemize}
\item[(a)] $|\mathcal{SS}(0,F)|= 8$;
\item[(b)] $\aleph_{\max}(0,F) = 2$ and $\aleph_{\min}(0,F) = 1$;
\item[(c)] $|\{X \in \mathcal{SS}(0,F): |\main(X)|=1\}|=4$ and $|\{X \in \mathcal{SS}(0,F): |\main(X)|=2\}|=4$;
\item[(d)] $\delta(0^+,F) = 0$ and $\Delta(0^+,F) = 4$;
\item[(e)] $\delta(0^-,F) = 0$ and $\Delta(0^-,F) = 1$;
\item[(f)] \begin{eqnarray*}
           \vert \{v \in V(F): d_{(0^+,F)}(v)=0\} \vert &=& 2, \\
           \vert \{v \in V(F): d_{(0^+,F)}(v)=1\} \vert &=& 0, \\
           \vert \{v \in V(F): d_{(0^+,F)}(v)=2\} \vert &=& 4, \\
           \vert \{v \in V(F): d_{(0^+,F)}(v)=3\} \vert &=& 0, \\
           \vert \{v \in V(F): d_{(0^+,F)}(v)=4\} \vert &=& 1;
           \end{eqnarray*}
\item[(g)] \begin{eqnarray*}
            \vert \{v \in V(F): d_{(0^-,G)}(v)=0\} \vert &=& 3,\\
            \vert \{v \in V(F): d_{(0^-,G)}(v)=1\} \vert &=& 4;
           \end{eqnarray*}
\item[(h)] Let $V^+_d$ and $V^-_{d}$ be, respectively, the subsets of vertices with $0$-main degree and $0$-non-main degree equal to $d$.
           \begin{enumerate}
           \item $V^+_0=\{f_5, f_6\}$ is an independent set;  $d_F(f_5)=d_F(f_6)=4$.
           \item $V^+_2=\{f_1, f_2, f_3, f_4\}$; the induced subgraph $F[V^+_2]$ is isomorphic to the cycle
                 $C_4$; $d_{F}(f_1)=d_{F}(f_4)=4$ and $d_F(f_2)=d_F(f_3)=2$.
           \item $V^+_4=\{f_7\}$;  $d_F(f_7)=2$.
           \item $V^-_0=\{f_5, f_6, f_7\}$; the induced subgraph $F[V^-_0]$ is isomorphic to the complete graph $K_3$; $d_F(f_5)=d_F(f_6)=4$ and $d_F(f_7)=2$.
           \item $V^-_1=V^+_2$.
           \end{enumerate}
\end{itemize}

\subsection{A comparison between $G$ and $F$}

The condition (h) of Proposition~\ref{properties_for_ip} fails to hold for $G$ and $F$
in the lists of vertex degrees obtained in items 1--4 of (h) and also for the induced subgraphs obtained in the item 4.

\bigskip

{\bf Acknowledgements} We are very grateful to the anonymous referees for their comments which improved the original manuscript. In particular, a referee suggested the third problem of Section~\ref{sec:problems}.

The second author is supported by the Center for Research and Development in Mathematics and Applications (CIDMA)
through the Portuguese Foundation for Science and Technology (FCT -- Fundação para a Ciência e a Tecnologia),
reference UIDB/04106/2020. The fourth author is supported by the Serbian Ministry of Education, Science and Technological Development via the University of Belgrade.

\end{document}